\documentclass[oneside,a4paper, 12pt]{article}

\usepackage{amsbsy, amssymb,  theorem}
\usepackage{pstricks,pst-plot,pst-node, hyperref}
\usepackage{epsf,newlfont,graphicx,amsbsy,stmaryrd}

\def\text{\hbox}

\usepackage{graphicx}
\usepackage{tikz}
\usepackage{caption}

\newtheorem{theorem}{Theorem}[section]
\newtheorem{lemma}[theorem]{Lemma}

\theoremstyle{definition}

\theoremstyle{remark}

\begin{document}

\title{Members of Thin $\Pi_1^0$ Classes and Generic Degrees \footnote{Stephan is supported by Singapore Ministry of Education Academic Research Fund Tier 2 grant MOE2016-T2-1-019 / R146-000-234-112.  Wu is supported by Singapore Ministry of Education Academic Research Fund Tier 2 grant MOE2016-T2-1-083 (M4020333); NTU Tier 1 grants RG32/16 (M4011672) and RG111/19 (M4012245). }}


\author{Frank Stephan\\ Department of Mathematics\\ Department of Computer Science\\ National University of Singapore\\ 10 Lower Kent Ridge, Singapore 119076\\
{\it email}: fstephan@comp.nus.edu.sg\\
\ $\ $\\
Guohua Wu\ \  Bowen Yuan\\
Division of Mathematical Sciences\\ School of Physical \&
Mathematical Sciences\\   Nanyang Technological
University\\ Singapore 637371\\
{\it email}: guohua@ntu.edu.sg;  yuan0058@e.ntu.edu.sg}

\maketitle

\begin{abstract}
	A $\Pi^{0}_{1}$ class $P$ is thin if every $\Pi^{0}_{1}$ subclass $Q$ of $P$ is the intersection of $P$ with some clopen set.
	In 1993, Cenzer, Downey, Jockusch and Shore initiated the study of Turing degrees of members of thin $\Pi^{0}_{1}$ classes, and proved that degrees containing no members of thin $\Pi^{0}_{1}$ classes can be recursively enumerable, and can be minimal degree below {\bf 0}$'$.
	In this paper, we work on this topic in terms of genericity, and prove that all 2-generic degrees contain no members of thin $\Pi^{0}_{1}$ classes. In contrast to this, we show that all 1-generic degrees below {\bf 0}$'$ contain members of  thin $\Pi^{0}_{1}$ classes.
\end{abstract}


\date{}


\section{Introduction}
In this paper, we will continue the study of Turing degrees of members of  thin $\Pi^{0}_{1}$ classes. Here a $\Pi^{0}_{1}$ class $P$ is thin if every $\Pi^{0}_{1}$ subclass $Q$ of $P$ is the intersection of $P$ with some clopen set.
Historically, thin classes were first constructed by Martin and Pour-El in their paper \cite{martin1970axiomatizable}, when they constructed an axiomatizable essentially undecidable theory such that any axiomatizable extension of it is a finite extension.
The concept of thin $\Pi^{0}_{1}$ class was first raised explicitly by Downey in his PhD thesis \cite{downeythesis}.
In a thin $\Pi^{0}_{1}$ class, computable elements are all isolated, and hence, perfect thin $\Pi^{0}_{1}$ classes contain no computable element.

In \cite{CDJS}, Cenzer, Downey, Jockusch and Shore considered the Cantor-Bendixson ranks of members of countable thin $\Pi^{0}_{1}$ classes, and then constructed an r.e. degree, and also a minimal degree below {\bf 0}$'$, whose elements are not members of any thin $\Pi^{0}_{1}$ class. We call such degrees thin-free degrees.
Downey, Wu and Yang recently proved in \cite{DWY2018} that the r.e. thin-free degrees are both dense and co-dense in r.e. degrees, providing another class of r.e. degrees which are both dense and co-dense in r.e. degrees. Recall that the first such a class is the set of of branching r.e. degrees, proved by Fejer \cite{Fejer} and Slaman \cite{Slaman}. Consequently, we call a Turing degree thin if it is not thin-free.

The construction of thin-free minimal degrees below {\bf 0}$'$  uses the $e$-splitting trees, which provides a framework for us to include and exclude infinitely paths of a given $\Pi_1^0$ class $P$, to show that $P$ is not a thin class. In the construction of a thin-free r.e. degree, they used an effective version of this idea, and constructed various intervals to capture the wanted infinitely paths.
The construction of thin-free minimal degrees below {\bf 0}$''$ turns out to be much easier as it is based on Spector's construction, instead of Sacks' construction. Based on this observation, in his  PhD thesis \cite{yuanthesis}, Yuan proved the existence of a hyperimmune-free minimal thin-free degree below {\bf 0}$''$.
We will use the same idea to show that nonrecursive sets below a 2-generic degree are thin-free. In particular, all 2-generic sets are of thin-free degree. In contrast to this, we show that all 1-generic degrees below {\bf 0}$'$ are not thin-free.

Our notation are standard. Most of the concepts and notation we use in paper can be found in books \cite{Cooperbook, Od, So}.

\section{2-generic sets are thin-free}

Recall that a set $A$ is $2$-generic if it meets or avoids every $\Sigma^{0}_{2}$ subset of $2^{<\omega}$.
Jockusch observed \cite{jockusch1980degrees}  that for any 2-generic set $A$, for any partial recursive functional $\Phi$ with $\Phi^A$  total and nonrecursive, $A$ has an initial segment $\sigma$ such that the set $T = \{ \rho: \rho \subseteq \Phi^{\tau}, \ \hbox{where}\ \tau\supseteq \sigma \}$ is a recursive extendible tree without isolated infinite paths.

\begin{lemma} [Jockusch]
	If $A$ is 2-generic and $\Phi^{A}$ is total and nonrecursive, then there is a $\sigma\subset A$ such that for all $\tau\supseteq\sigma$,
	\begin{itemize}
		\item for any $x$, $\tau$ has an extension  $\rho\supseteq\tau$ with $\Phi^{\rho}(x)\downarrow$,
		\item there is a $\Phi$-splitting extension above $\tau$.
	\end{itemize}
\end{lemma}

This is a property we need to show that a given $\Pi_1^0$ class is not thin.

\begin{theorem}
	Any nonrecursive set Turing below a 2-generic set is thin-free. Thus, 2-generic sets are of thin-free degree.
\end{theorem}

\medskip
\noindent {\bf Proof:}
	Let $A$ be a 2-generic set, and assume that $\Phi^{A}$ is total and nonrecursive. Suppose that $\Phi^{A}$ lies on a primitive recursive tree $P$, where $\Phi$ is a $\{0,1\}$-valued partial recursive functional. We will show that $[P]$ is not thin.
	
	Let $V$ be the set of strings $\tau$ such that:
	\begin{itemize}
		\item $\exists x \forall \rho\supseteq\tau[\Phi^{\rho}(x)\uparrow]$, or
		
		\item $\forall \rho, \pi\supseteq\tau\forall x[\Phi^{\rho}(x)\downarrow\ \& \ \Phi^{\pi}(x)\downarrow\ \Rightarrow\ \Phi^{\rho}(x) = \Phi^{\pi}(x)]$, or
		
		\item $\Phi^{\tau}\notin P$.
	\end{itemize}
	$V$ is $\Sigma^{0}_{2}$, and hence, $A$ either meets or avoids $V$, as $A$ is 2-generic.
	Note that $A$ cannot meet $V$, because $\Phi^{A}$ is assumed to be total, nonrecursive, and lies on $P$.
	So $A$ avoids $V$. This implies the existence of a string $\sigma\subset A$ such that for any $\tau\supseteq\sigma$,
	\begin{itemize}
		\item [(1)] for any $x$, $\tau$ has an extension  $ \rho\supseteq\tau$ with $\Phi^{\rho}(x)\downarrow$, and
		
		\item [(2)]  $\tau$ has two extensions $\rho, \pi $   such that for some $x$, $\Phi^{\rho}(x)\!\downarrow$, $\Phi^{\pi}(x)\!\downarrow$, and $\Phi^{\rho}(x) \neq \Phi^{\pi}(x)$, and
		
		\item [(3)] $\Phi^{\tau}\in P$.
	\end{itemize}
	
	For $\sigma$ above, consider  the tree $T = \{ \rho\subseteq \Phi^{\tau}:  \tau \supseteq \sigma\}$. We show that $T$ is extendible as follows.
	Suppose that for $\tau \supset \sigma$, $\Phi^{\tau}$ is a finite string on $T$. Then for any $x\geq|\Phi^{\tau}|$ (meaning that $\Phi^{\tau}(x)\uparrow$), by (1), there is some $\pi\supseteq\tau$ such that $\Phi^{\pi}(x)\downarrow$, which means that $\Phi^{\pi}$ extends $\Phi^{\tau}$. As $\Phi^{\pi}$ is on $T$,  $\Phi^{\tau}$ is extendible on $T$.
	
	Moreover, $T$ is recursive. Fix $l$. If $l \leq |\Phi^{\sigma}|$,  then only one string of length $l$, i.e., $\Phi^{\sigma}\upharpoonright l$, is on $T$.
	If $l > |\Phi^{\sigma}|$, let $\tau$ be the first string extending $\sigma$ such that $|\Phi^{\tau}| \geq l$. Such a $\tau$ exists by (1) and can be found recursively by enumerating strings extending $\sigma$.	Let $u(l) = \max\{ \varphi(x): x < l \}$, where $\varphi(x)$ is the use of $\Phi^{\tau}(x)$. $u(l)$ is recursive. By enumerating all strings of length $u(l)$ extending $\sigma$, we get all strings on $T$ of length smaller than $l$.
	
	Now consider the leftmost path $C$ though $T$.
	Since $T$ is a recursive, extendible tree, $C$ is recursive.
	$C$ lies on $T$, so $C \supset \Phi^{\sigma}$.
	Let $\rho_{0} = \sigma$, and for a given $\rho_{i}$, let $\rho_{i+1}$, $\tau_{i+1}$ and $x_{i+1}$ be the first triple $(\rho, \tau, x)$ such that $\rho,\tau\supseteq\rho_{i}$, $\Phi^{\rho}(x)\downarrow$, $\Phi^{\tau}(x)\downarrow$, $\Phi^{\rho}(x) \neq \Phi^{\tau}(x)$, and $\Phi^{\rho}\upharpoonright x \subset C\upharpoonright x$.	Such a triple exists because of (2) and the choice of $C$. Furthermore, as $C$ is recursive, the list of $\tau_{i}$ for $i \geq 1$ is recursive.
	By (3), $T$ is a subtree of $P$, so for any $i \geq 1$, there is an infinite path through $P$ extending $\Phi^{\tau_{i}}$.
	Let $S$ be the collection of all initial segments of $C$ and $\Phi^{\tau_{i}}$ for all $i \geq 1$, and the strings on $P$ extending $\Phi^{\tau_{i}}$ for  $i$ even.
	Then $S$ is a recursive subtree of $P$, which is not the intersection of $[P]$ with any clopen set.
	
	\begin{figure}[h]
		\centering
		\begin{tikzpicture}
		\filldraw (-0.5,1) circle (0.7pt);
		\filldraw (-1,2) circle (0.7pt);
		\filldraw (-1.5,3) circle (0.7pt);
		\filldraw (-1,3.5) circle (0.4pt);
		\filldraw (-1,3.7) circle (0.4pt);
		\filldraw (-1,3.9) circle (0.4pt);
		
		\coordinate [label=above:{$C$}] (a) at (-2,4);
		\coordinate [label=left:{$\Phi^\sigma$}] (a) at (-0.5,1);
		\coordinate [label=left:{$\Phi^{\rho_1}$}] (a) at (-1,2);
		\coordinate [label=right:{$\Phi^{\tau_1}$}] (a) at (-0.5,2);
		\coordinate [label=left:{$\Phi^{\rho_2}$}] (a) at (-1.5,3);
		\coordinate [label=right:{$\Phi^{\tau_2}$}] (a) at (-1,3);
		
		\draw (0,0) -- (-2,4);
		\draw [red][dashed] (-0.75,1.5) -- (-0.5,2);
		\draw [cyan] (-1.25,2.5) -- (-1,3);
		\end{tikzpicture}
		\captionof{figure}{The construction of $S$.}
	\end{figure}
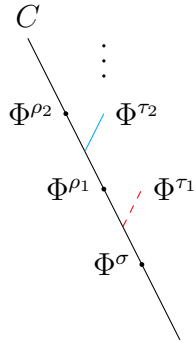
	
	So for any 2-generic set $A$, for all partial recursive functionals $\Phi$, if $\Phi^{A}$ is total, not recursive and $\Phi^{A} \in [P]$, then $[P]$ is not thin, completing the proof. \hfill $\Box$


\section{1-generic sets below {\bf 0}$'$ are not thin-free}

In this section, we show that all 1-generic degrees  below {\bf 0}$'$ are thin.
Let $A$ be a 1-generic set reducible to {\bf 0}$'$. To show a set $A$ is not thin-free, it suffices to construct a recursive thin tree $T$ containing a path $C$ which is Turing equivalent to $A$.

By Shoenfield's limit lemma, $A$ admits a $\Delta_2^0$ approximation. As $A$ is also 1-generic, $A$  admits a $\Delta_2^0$ approximations with an extra property, the so-called $\Sigma_1$-correctness. Here a recursive approximation $\{ \sigma_{s}: s\in\omega \}$ of $A$ is $\Sigma_{1}$-correct if for any infinite r.e. set $S$ of natural numbers, there exists some $s\in S$ such that $\sigma_{s}\subset A$.
Haught pointed out in \cite{Haught86} that this property is due to Shore, and used it to show that 1-generic degrees below {\bf 0}$'$ are downwards closed. For completeness of the paper, we present a proof of Shore's Lemma.

\begin{lemma} [Shore]
Any recursive approximation of a 1-generic set $A <_{T}${\bf 0}$'$ has a $\Sigma_{1}$-correct approximation.
\end{lemma}

\medskip
\noindent {\bf Proof:}
	As $A<_{T}${\bf 0}$'$, by Shoenfield's limit lemma, we can have a recursive appro\-xi\-ma\-tion of $A$, $\{ \sigma_{s} \}$ say.
	For an infinite r.e. set $S\subseteq\omega$, we define a set of strings $V = \{ \sigma_{s}: s\in S \}$, which is r.e..
	By 1-genericity, $A$ either meets or avoids $V$.
	$A$ cannot avoid $V$ because for any initial segment $\sigma$ of $A$, as $\{ \sigma_{s} \}$ approximates $A$, there exists some $s$ such that for all $t>s$, $\sigma_{t}\supseteq\sigma$. Thus, $A$ meets $V$, i.e., there is a $\sigma\subset A$ such that $\sigma \in V$, which implies the existence of $s\in S$ with $\sigma_{s}\subset A$.
\hfill $\Box$

\bigskip
In addition, Haught also pointed out that $A$ can actually have a $\Sigma_{1}$-correct approximation $\{ \sigma_{s}\}$  satisfying that $|\sigma_{s+1}| > |\sigma_{s}|$ for all $s\in\omega$. To see this, for any $\Sigma_{1}$-correct approximation $\{ \alpha_{s} \}$ of $A$, we define a function $f: \omega \rightarrow \omega$ inductively by taking $f(0) = 0$ and $f(s+1) = \mu t > f(s) (|\alpha_{t}| > |\alpha_{f(s)}|)$. Such a $t$ exists because $\{ \alpha_{s} \}$ is a recursive approximation. Note that $f$ is recursive and increasing.
Let $\sigma_{s} = \alpha_{f(s)}$, then $\{ \sigma_{s} \}$ is also a recursive approximation of $A$.
For any infinite r.e. set $S$, $V = \{ f(s): s\in S \}$ is r.e. and infinite. Since $\{ \alpha_{s} \}$ is $\Sigma_{1}$-correct, there is some $f(s)\in V$ such that $\alpha_{f(s)} \subset A$, which implies that $\sigma_{s} = \alpha_{f(s)} \subset A$. Thus $\{ \sigma_{s} \}$ is also $\Sigma_{1}$-correct.

\begin{theorem}
	A 1-generic degree {\bf a}$ <$ {\bf 0}$'$ is not thin-free.
\end{theorem}

\medskip
\noindent {\bf Proof:}
	Let $\{ \sigma_{s}: s\in\omega \}$ be a $\Sigma_{1}$-correct approximation of $A$ such that $|\sigma_{s+1}| > |\sigma_{s}|$ for each $s$. We will construct a recursive tree $T$, such that $[T]$ is thin and there is a path $C$ in $[T]$ with $C \equiv_{T} A$.

	Before we provide the construction of $T$, we first consider the set $S$ of all initial segments of $\sigma_{s}$, $s\in\omega$.
	$S$ is a tree since it is closed under initial segment, and $A$ is an infinite path through $S$. It is clear that $A$ is the only path on $S$.
	On the other hand, $S$ is r.e., as a string $\tau\in S$ if and only if   $\sigma_{s} \supseteq \tau$ for some $s$. As $\{ \sigma_{s}: s\in\omega \}$ is a $\Delta_2^0$ approximation of $A$, some strings on $S$ may not be extendible.
	
	We want the tree $T$ we are constructing to be recursive and extendible, and the construction of $T$ ``follows" the enumeration of $S$.
	With this in mind, we need an extra symbol $B$, standing for ``blank", such that $T$ is a subtree of $\{0,1,B\}^{<\omega}$ and all strings up to some length, $l(s)$ say, are defined on $T$ at each stage $s$. Here $l$ is a recursive function. It is clear that $T$ defined in this way can be coded into a binary tree recursively.
	
	For a finite string $\tau \in \{0,1,B\}^{<\omega}$ (or an infinite sequence $C \in \{0,1,B\}^{\omega}$, respectively), we let $\tau^{d}$ (or $C^{d}$) denote the string (or a finite string or an infinite subsequence, respectively) obtained by deleting all $B$ from $\tau$ (or $C$) while keeping the appearance of 0's and 1's the same order. For example, $(0010BB10B01)^{d} = 00101001$.
	
	\medskip
	\noindent {\bf Construction of $T$}:
	\begin{description}
		\item[ {\it Stage} $0$] Let $\emptyset$ be the root of $T$, and $l(0) = 0$.
		
		\item[ {\it Stage} $s+1$] $l(s)$, and all strings $\tau$ on $T$ of length $l(s)$ are already defined by the end of stage $s$.
		For $\rho$, a string on $T$ of length $l(s)$ with $\rho^{d} \subseteq \sigma_{s}$, let $l(s+1) = l(s) + m + 1$, where $m = |\sigma_{s}| - |\rho^{d}|$. Now there are $m+1$ steps for $i = 0,1,...,m$, and at step $i$, for strings $\tau$ on $T$ of length $l(s) + i$, if $\tau^{d} \subseteq \sigma_{s}$, put $\tau^{\smallfrown}0$ and $\tau^{\smallfrown}1$ into $T$, else put $\tau^{\smallfrown}B$ into $T$.
	\end{description}
	
	Note that for all $\tau$'s above, $\tau^{d}\upharpoonright |\tau^{d}|-1$ are already on $S$ by stage $s$, while none of $\tau^{d}$ is on $S$ yet, and	at stage $s+1$, with $\sigma_{s}$ just being put on $S$, for the string $\rho$ of length $l(s)$ on $T$ and $\rho^{d} \subseteq \sigma_{s}$, let $\pi$ be the string such that $\sigma_{s} = \rho^{d} \pi$ and $m = |\pi|$.
	Then (1) we put, if $m \geq 1$, $\rho^{\smallfrown}(\pi\upharpoonright i)^{\smallfrown}0$ and $\rho^{\smallfrown}(\pi\upharpoonright i)^{\smallfrown}1$ for $0 \leq i \leq m-1$, and $\rho^{\smallfrown}\pi^{\smallfrown}0$ and $\rho^{\smallfrown}\pi^{\smallfrown}1$ into $T$ in order. In this manner, all strings ending with $0$ or $1$ on $T$ longer than $l(s)$ are defined.
	(2) we extend all other strings on $T$ by $B$ up to length $l(s) + m + 1$. This completes the construction of $T$ up to length $l(s) + m + 1$.
		
	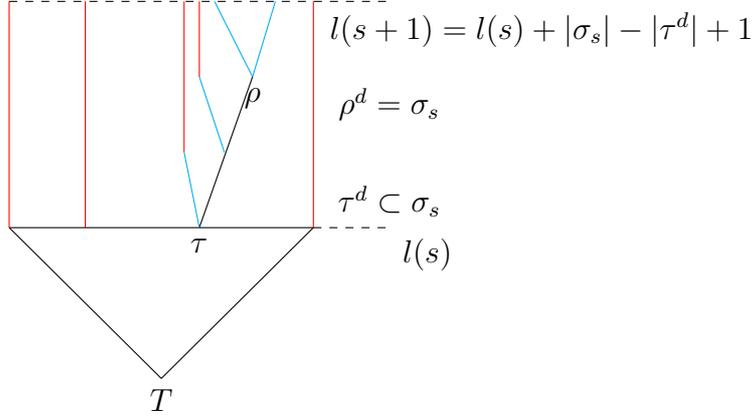
\begin{figure}[h]
		\centering
		\begin{tikzpicture}	
		\coordinate [label=below:{$T$}] (a) at (0,0);
		\coordinate [label=below:{$\tau$}] (a) at (0.5,2);
		\coordinate [label=below:{$\rho$}] (a) at (1.2,4);
		\coordinate [label=below:{$l(s+1)= l(s) + |\sigma_{s}| - |\tau^{d}| + 1$}] (a) at (5,5);
		\coordinate [label=below:{$l(s)$}] (a) at (3.5,2);
		\coordinate [label=above:{$\tau^d\subset\sigma_s$}] (a) at (3,2);
		\coordinate [label=below:{$\rho^d=\sigma_s$}] (a) at (3,4);
		
		\draw (0,0) -- (-2,2);
		\draw (0,0) -- (2,2);
		\draw (-2,2) -- (2,2);
		\draw [red](-2,2) -- (-2,5);
		\draw [red](-1,2) -- (-1,5);
		\draw [red](2,2) -- (2,5);
		\draw (0.5,2) -- (1.2,4);
		\draw [cyan] (1.2,4) -- (0.7,5);
		\draw [cyan] (1.2,4) -- (1.5,5);
		\draw [cyan] (0.5,2) -- (0.3,3);
		\draw [cyan] (0.83,3) -- (0.5,4);
		\draw [red] (0.3,3) -- (0.3,5);
		\draw [red] (0.5,4) -- (0.5,5);
		\draw [dashed] (-2,5) -- (3,5);
		\draw [dashed] (2,2) -- (3,2);
		\end{tikzpicture}
		\captionof{figure}{An illustration of construction of $T$.}
	\end{figure}
	
	What does $T$ look like?
	Consider all the strings on $T$ of length $s+1$, i.e., those strings on the $(s+1)$-st level of $T$.
	There are exactly $s+2$ strings on the $(s+1)$-st level, and among them, two strings end with $0$ and $1$ and share the common initial segment of length $s$, and the other strings end with $B$.
	Let $\tau_{i}$ for $0\leq i\leq s+1$ be the strings on the $(s+1)$-st level of $T$, then $\tau_{i}^{d}$ are incompatible with each other, and only one of them, say $\tau_{j}^{d}$, is extendible on $S$ (equivalently, $\tau_{j}^{d} \subset A$). So above the $(s+1)$-st level of $T$, there are only finitely many paths above $\tau_{i}$ on $T$ for each $i \neq j$, and infinitely many paths on $T$ extending $\tau_{j}$.
	Since $A$ is the unique path on $S$, we know that there is only one path $C$ on $T$ which contains infinitely many $0$ and $1$, with $(C)^{d} = A$. For other paths $D$, there is some $n$ such that for all $x>n$, $D(x) = B$, and $(D)^{d}$ is a string not in $S$.
	
	\begin{figure}[h]
		\centering
		\begin{tikzpicture}	
		\coordinate [label=below:{$T$}] (a) at (0,0);
		\coordinate [label=above:{$C$}] (a) at (0.5,4);
		
		\draw (0,0) -- (-1,1);
		\draw [red] (0,0) -- (1,1);
		\draw (-1,1) -- (-1,2);
		\draw [red] (1,1) -- (0.5,2);
		\draw (1,1) -- (1.5,2);
		\draw [dashed](-1,2) -- (-1.2,2.2);
		\draw [dashed](-1,2) -- (-0.8,2.2);
		\draw [dashed](-0.8,2.2) -- (-0.8,4);
		\draw [dashed](-1.2,2.2) -- (-1.2,4);
		\draw [dashed](1.5,2) -- (1.7,2.2);
		\draw [dashed](1.5,2) -- (1.3,2.2);
		\draw [dashed](1.7,2.2) -- (1.7,4);
		\draw [dashed](1.3,2.2) -- (1.3,4);
		\draw [dashed] (0.5,2) -- (0.1,4);
		\draw [dashed] (0.5,2) -- (0.9,4);
		\draw [red] (0.5,2) -- (0.5,4);
		\end{tikzpicture}
		\captionof{figure}{An example of how $T$ and $C$.}
	\end{figure}
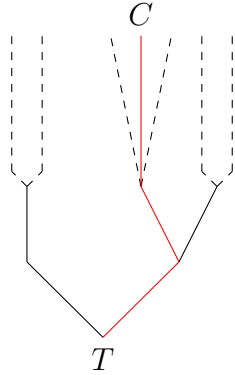
	
	\medskip
	$A \leq_{T} C$ by $(C)^{d} = A$.
	$C \leq_{T} A$ because for any $x$, there is some $s_{x}$ such that $l(s_{x}) \leq x \leq l(s_{x}+1)$.	Then find the least $s>s_{x}$ such that $\sigma_{s} \subset A$. Such an $s$ exists since $\{\sigma_{s}\}$ is $\Sigma_{1}$-correct.
	For this $s$, find $\tau$ on $T$ of length $l(s)$ with $\tau^{d} \subset A$. Then $\tau$ is an initial segment of $C$, and $\tau(x) = C(x)$.
	
	We now show that $T$ has a neat ``splitting'' property.
	
	\begin{lemma}
		For any $n$, there is some $m>n$, such that for strings on the $n$-th level of $T$, all but one string do not split above level $m$ on $T$.
	\end{lemma}

\medskip
\noindent {\bf Proof:}
		For any $n$, let $\tau_{i}$ for $0\leq i \leq n$ be the strings on the $n$-th level of $T$, and $k = \max_{i}\{|\tau_{i}^{d}|\}$.
		For this $k$, by $\Sigma_{1}$-correctness, there is a stage $t \geq k$ such that $\sigma_{t}\subset A$. Then there is some stage $s$ such that for all stages $s'>s$, $\sigma_{s'}\supset\sigma_{t}$.
		
		Let $m = l(s+1) - 1$ and $\rho$ be the string on the $m$-th level of $T$ such that $\rho^{d} = \sigma_{s}$.
		Then $\rho^{d} \supset \sigma_{t}$, and after stage $s$, all strings on the $m$-th level of $T$ except $\rho$ can only be extended by $B$, and thus if $\tau_{i} \nsubseteq \rho$, $\tau_{i}$ does not split above level $m$ on $T$.\hfill $\diamond$
	
\bigskip
	Thus, all strings on $T$ are extendible, and $[T]$ contains infinitely many paths.
	
	The $\Sigma_{1}$-correctness of $\{\sigma_{s}\}$ guarantees that $T$ is thin. Let $U$ be a recursive subtree of $T$. There are two cases.
	
	\begin{enumerate}
		\item[(1)] $C\not\in [U]$. If so, then there is some $n$ such that $C\upharpoonright n \notin U$. As $C^{d} = A$, $(C\upharpoonright n)^{d} \subset A$, and thus there are only finitely many paths on $U$. Suppose that up to length $m$, strings on $U$ do not split, let $N$ be the union of cones above the strings on $U$ of length $m$, then $N$ is clopen and $[U]$ is the intersection of $[T]$ with $N$.
		
		\item[(2)] $C\in [U]$. Consider $V = \{ s: \tau\in T\backslash U, \tau^{d} \in S \text{ for some stage }\ s\}$. Then by the $\Sigma_{1}$-correctness,  $V$ is finite. Otherwise, there is some $s\in V$ such that $\sigma_{s} \subset A$, which means that there is some $\tau$ in $T\backslash U$ such that $\tau^{d}=\sigma_{s}$ is enumerated into $S$ at stage $s$. By the construction, $\tau^{d} \subset A$, so $\tau \subset C$. However, $C\in[U]$ implies that for all $n \geq 0$, $C\upharpoonright n \in U$, a contradiction.
		
		Since $V$ is finite, there are finitely many paths on $T$ but not on $U$. Thus, there is a clopen set $N$ such that $[T]\backslash[U] = [T] \cap N$, and hence $[U]=[T]\backslash N=[T]\cap \overline{N}$. The complement of $N$ is what we need.	\end{enumerate}

	\medskip
	This shows that among these infinitely paths, exact one path is Turing equivalent to $A$.
	Thus, we obtain a path $C$ in a thin class $[T]$ with $C \equiv_{T} A$, and $A$ is not of thin-free degree.	\hfill $\Box$


\bibliographystyle{amsplain}

\end{document}